# NEW ENTROPY AND OTHER STATE FUNCTIONS FROM ANALYSIS OF OPEN SYSTEM CARNOT CYCLE

Christopher Gunaseelan Jesudason
Chemistry Department, University of Malaya
50603 Kuala Lumpur, Malaysia
Email: jesu@um.edu.my, Fax:03-79674193

**ABSTRACT**
A first principles analysis of an open system thermodynamical Carnot cycle is provided, and the results are compared to those proposed by Gibbs for open systems. The Kelvin-Clausius statement concerning heat transfer for reversible cycles is taken as an axiom, from which several rigorous theorems are proven. An equation is derived that resembles a Gibbs-Duhem relation relating convected entropies, from which two distinguishable forms of entropy are proven to exist for such systems, which questions prevailing developments which presume a singular or characteristic entropic form which couple all work and heat flows, such as in Onsager first order thermodynamics. In particular, a closed path undergone by the system does not return the environment to the initial state for one of these entropic forms. The Biot assertion that the entropy contribution due to diathermal heat transfer does not form a state function is therefore contradicted and a local entropy is shown to exist. Several other new composite state functions for work and heat flow are shown to exist for open systems. From Gibbs' results, it is suggested that the ensuing chemical potentials used routinely may possibly ignore heat effects. The functions developed here are suitable for application since the functions are proven to exist, rather than presumed to exist.

## 1. Introduction
Gibbs [1] did not carry out a rigorous open system Carnot engine analysis; he presumed that one could extend the entropy concept through the presumed existence of the function of state without proving existence. More recently, [2] Biot has attempted to extend this concept to open systems by using the concept of a primary cell and supply cells with a singular Thermal Well (TW) or reservoir, but the concept is not local in nature, since the variables, although pertaining to the primary cell, refer to extraneous components. He also used the traditional notion of heat pumps from closed Carnot engines. Here, we analyze an open Carnot engine, where the TW is replaced by the heat reservoirs at different temperatures as in the original Carnot elaborations having a local characteristic, rather than one extended arbitrarily into space as in Biot's method. The result of this analysis contradicts Biot rigorous analysis in at least one major way concerning conductive heat transport, and some new functions of state are deduced, of immediate significance in both equilibrium and nonequilibrium studies.

## 2. System Description

The elementary open Carnot engine is described in Fig. 1 below where $ab, cd$, are segments of an isothermal path and $bc, da$ are adiabats, and if masses $\delta m_i$ are injected and extracted reversibly along an isothermal segment $a'b$ and $c'd$ respectively through a selectively permeable membrane, the cycle is termed $C_{iso}$ whereas if the same process is carried out adiabatically along adiabatic segments $bb', dd'$, the cycle is termed $C_{adia}$. The system itself may or may not have chemical reactions occurring within, but possesses a minimum independent set of variables $\Sigma$ defined $\Sigma = \{\mathbf{P}, \mathbf{V}, \mathbf{m}\}$ where the listed variables denote the system intensive and extensive variables and the mass amounts present respectively, where in particular $\mathbf{m} = \{m_1, m_2...m_i...m_n\}$ is the set of 'pure' substances injected into the open system and which might consequently react; $\mathbf{m} \cup \mathbf{V}$ are the extensive quantities such that $-\mathbf{P} \cdot (\delta \mathbf{V} + \delta \mathbf{m})$ represents the work gained by the system for arbitrary displacements $(\delta \mathbf{V} + \delta \mathbf{m})$. Thermodynamical points $a, b, c, d$ also represents the state variables $\Sigma$ for those labeled states. For $C_{iso}$, masses are injected reversibly at constant temperature by the simultaneous exchange of heat via the system boundaries, whereas for $C_{adia}$ there is no temperature control for the system during reversible exchange of masses. The temperatures along isotherms $ab$ and $cd$ is $T_1$ and $T_2$ respectively; and any of these cycles, heat is exchanged with the thermal reservoirs labeled $T_1$ and $T_2$ held respectively at the temperatures $T_1$ and $T_2$ (here the

s

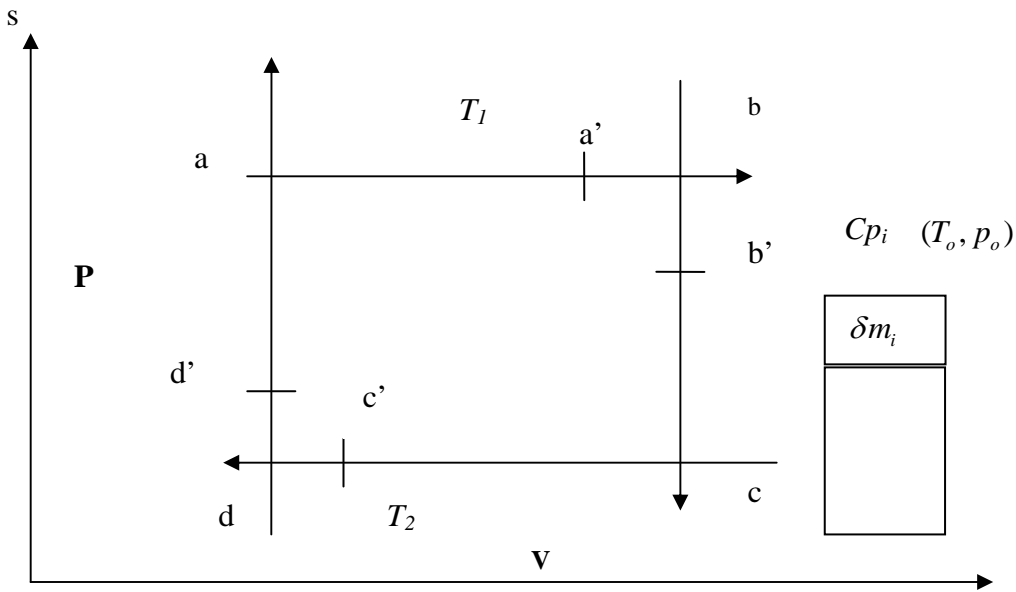

**Figure 1: A system taken through a reversible cycle along abcda in multidimensional space $(\mathbf{P}, \mathbf{V}, \mathbf{m})$ with mass exchange $\delta m_i$.**

symbol refers to both variable and item, depending on context). Similarly $C_{iso}$ and $C_{adia}$ also refer to the systems undergoing the particular type of transition. The supply cells $Cp_i$ for substance $i$ are all held at a standard temperature and pressure denoted $ss$ at $(T_o, p_o)$. Without loss of generality, we assume $T_1 > T_2 > T_0$. In some cases pure substances cannot exist in the standard state, for instance in the dimer reaction $2A \rightleftharpoons A_2$ where both species $A$ and $A_2$ are present in principle for all thermodynamical states, from the classical expression $RT \ln K = -\Delta_r G^{\ominus}$ where $K$ is the equilibrium constant, $\Delta_r G^{\ominus}$ the standard free energy change, $R$ the gas constant and $T$ the temperature. In such instances, the standard state would be one where *both* species coexisted in the supply cell, where if a semi permeable membrane were used to inject species $A$ into the system, then the reversible work done in this case would be for the combined masses such that $\delta m_k = \delta m_{A,k} + \delta m_{A_2,k}$ where the incremental work term would involved both the injection of mass against the partial pressure $p_k$ for species $A$, but also the work done for the separation of the molecules, which involves the isothermal absorption of heat through the system. The two elementary cycles $C_{iso}$ and $C_{adia}$ and their combinations will be examined in turn. The term 'heat absorbed' or 'work' refer to the work done on the system and the heat absorbed by the system only unless stated otherwise. This system termed $S$ has portions of its surface $\partial S_{iso}$ that is diathermal, and other portions $\partial S_{adia}$ which can exchange particles through a selectively permeable membrane, on which surface work energy may be expended through reversible pressure and surface tension forces; the particles exchanged are in the pure state (or mixed state in reactive pure substances) at the opposite side of the membrane away from the system, and the walls of this device (of vanishingly small heat capacity) may transfer and exchange heat with $\partial S_{iso}$ through thermal conduction. The entire hypersystem which consists of $S \cup T_2 \cup T_1 \cup \left[ \bigcup_{i=1}^{n} Cp_i \right]$ exists in an inertial frame where the environment is stationary, and where the work and heat exchange with environment and subsystems is such that no relative motions are allowed that would affect the energy interactions described. Without loss of generality, we may assume that the fluids are compressible as $\delta m_i \to 0$.

**3. Relations Amongst the Variables**
(i) $C_{iso}$ Cycle Analysis
The reversible heat and work terms undergone by this cycle are as follows:
$\delta Q^{T_1}_{a-a',sys}$ is the heat absorbed along isotherm $aa'$ and similarly the heat absorbed along isotherm $cc'$ is $-\delta Q^{T_2}_{c-c',sys}$. The energy per unit mass of substance $i$ in the supply cell is denoted $U^o_{ss,i}$. The partial pressure of component $i$ in the system which is at equilibrium with a pure component at the same pressure is a function of variables $\Sigma$, $p_i = p_i(\Sigma)$ and the work done per unit mass of injecting $i$ is $\zeta_{iso,i}(\Sigma)$ at constant $\{\mathbf{P}, m_i \notin \mathbf{m}\}$ and the work of injecting amount $\delta m_i$ is

$$\delta W_{inj,i} = \zeta_{iso,i}(\Sigma) dm_i \quad (1)$$

and the reversible heat absorbed per unit mass at the indicated temperature $T$ is $\phi_i(\Sigma)$ and for amount $\delta m_i$ is

$$\delta Q^T_{inj,i} = \phi_i(\Sigma) dm_i \quad (2)$$

Then at points between $c'd'$ the work of extraction $\delta W_{ext,iso}$ is

$$\delta W_{ext,iso,i}(\Sigma) = -\delta W_{inj,iso,i}(\Sigma) = -\zeta_{iso,i}(\Sigma_{c'd'}) \quad (c' \to d') \quad (3)$$

and is accompanied by heat of extraction (released to $T_2$)

$$\delta Q^{T_2}_{ext,i} = -\phi_i(\Sigma_{c'd'}) dm_i = -\delta Q^{T_2}_{inj,i}. \quad (4)$$

Similarly, the work done per unit mass for the adiabatic compression at any point $\Sigma$ without heat exchange is $\zeta_{adia,i}(\Sigma)$ and

$$\delta W_{ext,i} = -\delta W_{inj,i} = -\zeta_{adia,i}(\Sigma)\delta m_i \quad (5)$$

where the appropriate subscripted form of $\delta W_{ext,i}$ is used depending on whether an adiabatic or isothermal mass transfer is implied. The work done by the environment on the system along $aa'$ and $cc'$ is respectively $\delta W^{T_1}_{aa'}$ and $\delta W^{T_2}_{cc'}$ and heat absorbed is $-\delta Q^{T_2}_{cc'}$ along $cc'$; the work done by the environment on system along $bc, da$ is respectively $\delta W_{bc}$ and $\delta W_{da}$; the transition being adiabatic, there is no reference to heat absorption or temperature. The work of extraction per unit mass from $Cp_i$ is $\zeta_{Cp,i}(p_o, T_o)$ and the work done by the environment is given by

$$\delta W_{ext,Cp,i} = \zeta_{Cp,i}(p_o, T_o)\delta m_i. \quad (6)$$

$W_{ss,i}$ refers to the total work per unit mass to extract $i$ from $Cp_i$ at state $ss$. Biot[2] over three decades developed a theory where the heat source for the heat pumps are given by a thermal well $T_w$ at $T_o$ for the whole hypersystem; we dispense with this $T_w$ and use the Carnot idea of a minimum of two thermal wells at $T_1$ and $T_2$. The heat given up by $T_1$ by a reversible transfer of heat by heat pumps from the standard state to the state $(p_i, T)$ at equilibrium with $S$ through a semi permeable membrane for substance $i$ is

$$\Delta Q_i^{(p_i,T)} = T\left(\int_{ss}^{(p_i,T)} \frac{dQ_i}{T}\right)\delta m_i \quad (7)$$
$$= T\delta m_i \Delta S_{i,ss}^{(T,p_i)}$$

where $dQ_i = Cp\, dv + C_v dp$ is the heat absorbed per unit mass in the pressure volume $(p,v)$ path from the standard state $ss$ to $(p_i,T)$, and $\Delta S_{i,ss}^{(T,p_i)} = \int_{ss}^{(p_i,T)} \frac{dQ}{T}$ is the convected entropy relative to state $ss$. In addition to the heat pumps, the increment of work done per unit mass on the element $i$ along the pathway from $ss$ to the final state $(p_i, T)$ is $\Phi_i$, an imperfect differential, so that the increment of work $dW_i^{form}$ along the pathway for element $\delta m_i$ can be written

$$dW_i^{form} = d\Phi_i \delta m_i = (\mu_{p,i} dT + \mu_{T,i} dp)\delta m_i \quad (8)$$

and integrating along the path from $ss$ to $(p_i, T)$ yields $\Delta W_i^{form}(p_i, T)$ where $\mu_{p,i}$ and $\mu_{T,i}$ do not allow $\Delta W_i^{form}$ to be an exact differential.

Biot[2] deduced the existence of a work state function $\mathcal{V}$ which he used extensively in his thermodynamics, which is the total work done on heat pumps and the reversible work on the primary cell and elements of the supply cell. This result is a trivial consequence of the Kelvin-Clausius statement of the Second law, since only one thermal well $T_w$ is considered; if for instance there exists two different values $\mathcal{V}_1$ and $\mathcal{V}_2$ for the same final state, driving one system against another appropriately for two duplicate systems initially at the above same state to another where all materials are returned to the supply cells would lead to the pure conversion of heat in $T_w$ to work in a cycle with no other effect. There is an analogous simple result to the above.

**Lemma 1:** The reversible heat absorbed function for any substance $i$ $\Delta Q_i^{(p_i,T)} = T \delta m_i \Delta \mathcal{S}_{i,ss}^{(T,p_i)}$ and the total reversible heat absorbed function $\Delta Q_{total}^{\Sigma} = \sum_{i=1}^{n} T m_i \Delta \mathcal{S}_{i,ss}^{(T,p_i)}$ are perfect differentials or state variable functions

*Proof*
Since $\Delta \mathcal{S}_{i,ss}^{(T,p_i)}$ which is not coupled to any **m** is a state function from the second law, and so are **m** and $T$, it follows that $\Delta Q_i^{(p_i,T)}$ is a function of state and by simple integration of **m** and summation of all species implies $\Delta Q_{total}^{\Sigma}$ is a state function ∎

$\Delta Q_{total}^{\Sigma}$ may be used as a state function and expanded as such to yield relations based on the Maxwellian reciprocity condition.

**Lemma 2**: There exists work functions of state $\Delta \mathcal{W}^{\Sigma}$ and $\Delta \mathcal{W}_i^{\Sigma}$
*Proof*
The work done on the environment for the pumping mechanism of species $i$ is per unit mass $\Delta \mathcal{V}_{i,ss}^{(T,p_i)}$ is such that

$$\Delta \mathcal{V}_{i,ss}^{(T,p_i)} \delta m_i = T \delta m_i \Delta \mathcal{S}_{i,ss}^{(T,p_i)} - \Delta Q_{i,ss}^{(T,p_i)} \delta m_i \quad (9)$$

where $\Delta Q_{i,ss}^{(T,p_i)} = \int_{ss}^{(T,p_i)} dQ_i$ is the total heat absorbed by the element per unit mass. The energy change per unit mass is $\Delta U_{i,ss}^{\Sigma}$, is a function of state where

$$\Delta U_{i,ss}^{\Sigma} \cdot \delta m_i = \left( \Delta W_i^{form}(p,T) + \Delta Q_{i,ss}^{(T,p_i)} \right) \cdot \delta m_i \quad (10)$$

so there exists a function of state of work energies $\Delta \mathcal{W}_i^\Sigma$ given by comparing (9) and (10) to yield
$$\Delta \mathcal{W}_i^\Sigma = \Delta W_i^{form}(p,T) - \Delta \mathcal{V}_{i,ss}^{(T,p_i)} \quad (11)$$
and a total system work energy state function given as
$$\Delta \mathcal{W}^\Sigma = \sum_{i=1}^{n} \Delta \mathcal{W}_i^\Sigma \quad . \quad (12)$$
This function can also yield thermodynamical relations through Maxwellian reciprocity ∎

In general, we define the variation $\delta X$ of the above quantities $\Delta X$ as $\delta X_i = \Delta X_i \delta m_i$.

For the system undergoing a complete cycle so that $Cp_i$ and $S$ are restored to their original condition, the total work done on the environment by the system $\Delta W_{tot,iso}$ and by the heat pumps, all of which work cyclically are as follows

$$\Delta W_{tot,iso} = -\delta W_{aa'}^{T_1} - \delta W_{cc'}^{T_2} - \delta W_{bc} - \delta W_{da} + \delta \mathcal{V}_{ss}^{(T_1,p_{i,\mathbf{b}})}$$
$$- \delta \mathcal{V}_{ss}^{(T_2,p_{i,\mathbf{d}})} + \Delta W_i^{form}(p_{i,\mathbf{b}},T_1) - \Delta W_i^{form}(p_{i,\mathbf{d}},T_2)$$
$$- \delta W_{inj}(a'b) + \delta W_{inj}(c'd).$$
(13)

Total heat *lost* at $T_1$ reservoir $\Delta Q_{1,tot,iso}$ is
$$\Delta Q_{1,tot,iso} = \delta Q_{aa',sys}^{T_1} + \delta Q_{inj,i}^{T_1} + T_1 \delta m_i \Delta \mathcal{S}_{i,ss}^{(T_1,p_i)}. \quad (14)$$

Total heat *gained* at $T_2$, $\Delta Q_{2,tot,iso}$ is

$$\Delta Q_{2,tot,iso} = \delta Q_{c'c,sys}^{T_1} - \delta Q_{inj,i}^{T_2} + T_2 \delta m_i \Delta \mathcal{S}_{i,ss}^{(T_2,p_{i,c'd})}$$
$$= \delta Q_{cc',sys}^{T_1} - \delta Q_{inj,i}^{T_2} + T_2 \delta m_i \Delta \mathcal{S}_{i,ss}^{(T_2,p_{i,c'd})}.$$
(15)

The states as given in Fig.1 refer to the states a to d with changes in the $\{\mathbf{P},\mathbf{V}\}$ and where in particular the **m** variables have the following values: a: $\{\mathbf{m}=\mathbf{m}'\}$, b: $\{\mathbf{m}=\mathbf{m}'$, but $m_i = m_i' + \delta m_i'\}$, c: $\{\mathbf{m}=\mathbf{m}'$, but $m_i = m_i' + \delta m_i'\}$, d: $\{\mathbf{m}=\mathbf{m}'\}$. In this complete cycle we traverse the loop $a \to b \to c \to d$. A simple result is required to complete the proof of function of state.

**Lemma 3**: A closed system $S'$, subjected to the Kelvin-Clausius Axiom and which undergoes a totally reversible cycle exchanging heat at reservoirs at temperatures $T_1$ and $T_2$ is a Carnot device.
*Proof*
Trivially, if $T_1 > T_2$ and $S'$ absorbed heat $Q_1'$ at $T_1$ and ejected heat $Q_2'$ at $T_2$ with work $W' = Q_1' - Q_2'$ done on the environment, then we could run a Carnot device such that it gave up heat $Q_1'$ at $T_1$ and absorbed heat $Q_2$ at $T_2$ with work $W$ done on the environment. Then if $Q_2' \neq Q_2$, $-W = Q_1' - Q_2$ or $W' \neq -W$, and so the net effect (by choosing $S'$ to run in the appropriate direction) would be the total conversion of heat energy from $T_2$ to work with no other effect, violating the Kelvin-Clausius postulate. So $S'$ is Carnot ∎

(ii) $C_{adia}$ Cycle Analysis

As $a' \to b, c' \to d$, the reversible cycle traversed $a \to b \to c \to d \to a$ leads to the total work and heat absorbed by both reservoirs (subscripted adia. refers to the $C_{adia}$ cycle and other subscripts refer to the path or state).

$$\Delta W_{tot,adia} = -\delta W_{ab}^{T} - \delta W_{cd}^{T_2} - \delta W_{b'c} - \delta W_{d'a} + \delta \mathcal{V}_{ss}^{(T_1, p_{i,\mathbf{bb'}})}$$
$$-\delta \mathcal{V}_{ss}^{(T_2, p_{i,\mathbf{d'}})} + \Delta W_i^{form}(p_{i,\mathbf{b}}, T_1) - \Delta W_i^{form}(p_{i,\mathbf{d}}, T_2).$$
(16)

$$\Delta Q_{1,tot,adia} = \delta Q_{ab,sys}^{T_1} + \delta Q_{inj,i}^{T_1} + T_1 \delta m_i \Delta \mathcal{S}_{i,ss}^{(T_1, p_i)}. \quad (17)$$

$$\Delta Q_{2,tot,adia} = \delta Q_{dc,sys}^{T_2} + \delta Q_{inj,i}^{T_1} + T_1 \delta m_i \Delta \mathcal{S}_{i,ss}^{(T_1, p_i(c'd))}. \quad (18)$$

We note that there can be no heat absorption about an adiabatic segment.

(iii) $C_{comb}$ Cycle Analysis

$C_{comb}$ refers to a system where mass is injected(extracted) at an isothermal pathway e.g. ab and extracted(injected) at the adiabatic pathway e.g. da. For this reversible cycle, we still derive $\Delta W_{tot,comb}, \Delta Q_{1,tot,comb}$ and $\Delta Q_{2,tot,comb}$, as above.

Deduction for set $\mathbf{q} = \{C_{iso}, C_{adia}, C_{comb}\}$

For any q in $\mathbf{q}$ above, it is clear from Lemma 3 that because the cycles are all reversible and have as net effect the exchange heat at $T_1, T_2$ only, these cyclical paths must be Carnot and hence from the results of the Carnot system, it follows that

$$\frac{\Delta Q_{1,tot,q}}{T_1} - \frac{\Delta Q_{2,tot,q}}{T_2} = 0 \quad (19)$$

or

$$\frac{\Delta Q_{1,tot,q}}{\Delta Q_{2,tot,q}} = \frac{T_1}{T_2} \quad (20)$$

We have shown that any set in $\mathbf{q}$ fulfills the Carnot theorem which is defined as the conditions of efficiency given by (20) when heat is exchanged between two ports in a cycle. The above however have input/output of mass elements $i$ in the vicinity of the corners of the Carnot cycle; to place these ports elsewhere is equivalent to concatenating them with closed systems with matching adiabatic or isothermal pathways.

For an isothermal input/output system $S \equiv C_{iso}$, the concatenation (Fig. 2, I) produces a new Carnot engine system $S2 \cup S \cup S1$ with isothermal input/output ports at C and G at arbitrary distance from adiabatic paths DE and AH, where BG and CE are matching adiabatic paths. Similarly, for adiabatic inputs at C and G ( Fig. 2, II) of the original system $S$, we can concatenate this system with S3 and S4, with common isothermal pathways HC and GD for the composite system $S \cup S3 \cup S4$ where the ports at C and G are at arbitrary distance from the isotherms AB and FE. The proof that these composites are Carnot follow.

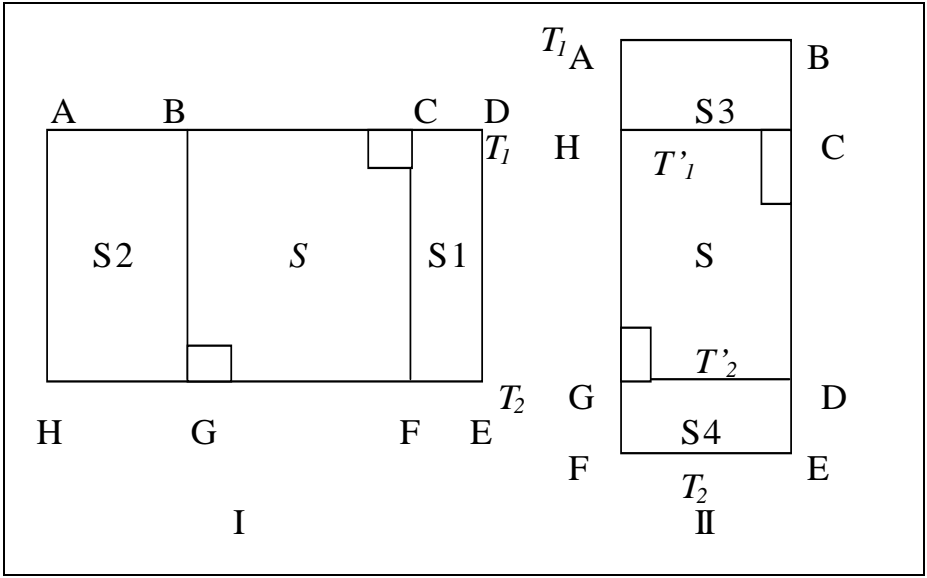

**Figure 2. Concatenating two closed systems S1 and S2 for system I. and S3 and S4 for system II to produce a new elementary Carnot engine $S_{new}$ with configuration $S2 \cup S \cup S1$ and $S \cup S3 \cup S4$ respectively, which have matching adiabatic boundaries (system I) or isothermal boundaries (system II) working between thermal reservoirs at the indicated temperatures.**

**Lemma 4**: The devices $S_{new}$ above are Carnot.
*Proof*
The proofs are basic. For $S_{new}$ (I), since S1, S and S2 are Carnot, then

$$\frac{Q_{1,S}}{T_1} + \frac{Q_{1,S2}}{T_1} + \frac{Q_{1,S3}}{T_1} = \frac{Q_{1,S_{new}(I)}}{T_1} = \frac{Q_{2,S}}{T_2} + \frac{Q_{2,S2}}{T_2} + \frac{Q_{2,S3}}{T_2} = \frac{Q_{2,S_{new}(I)}}{T_2}.$$ Since

$\frac{Q_{1,S_{new}(I)}}{T_1} = \frac{Q_{2,S_{new}(I)}}{T_2}$, the composite system is Carnot. For $S_{new}$ (II), since S3, S4, and

S are Carnot, then $\frac{\delta Q_{AB}}{T_1} = \frac{\delta Q_{HC}}{T'_1} = \frac{\delta Q_{GD}}{T'_2} = \frac{\delta Q_{FE}}{T_2}$. Since $\frac{\delta Q_{AB}}{T_1} = \frac{\delta Q_{FE}}{T_2}$ for

$S_{new}$ (II), it is Carnot. Similarly it can be shown that $C_{comb}$ for arbitrarily positioned input/output ports are Carnot∎

The above pertain to elementary finite cycles. We extend the above results to general cycles.

**4. System Entropies**

**(i)** Derivation of composite entropy

**Theorem 1.**: There exists a perfect differential $d\mathcal{S} = \dfrac{dQ_{tot}}{T}$ for the entropy function $\mathcal{S} = \mathcal{S}(\mathbf{P}, \mathbf{V}, \mathbf{m})$.

*Proof*

We have described the most general elementary Carnot cycles above. Consider an isothermal path connecting state $\mathbf{a} = (\mathbf{P}, \mathbf{V}, \mathbf{m})$ to $\mathbf{b} = (\mathbf{P'}, \mathbf{V'}, \mathbf{m'})$ such that state points $1, 2 ... N' = \mathbf{b}$ lie on the path where each state point $j$ is indexed by $j$ and where none of the points are coincident. The segment $j$ are the set of points on the path bounded by the point $j+1$ and $j$ less the state point $j$. The reversible heat given up by a reservoir at temperature $T_{1,j}$ along segment $j$ is denoted $\delta Q_{tot,j}$, and the heat absorbed by the thermal reservoirs along segments $j'$ is denoted $\delta Q_{tot,j'}$ at temperature $T_{2,j'}$. Points $j, j+1, j'+1, j'$ are the The function boundary points for an elementary Carnot cycle $S_j$ of species $C_{comb}$ as previously described. Also the points $j$ are placed in increasing order with a distance metric from the point $\mathbf{a}$ and $j'$ is a point connected to $j$ along an adiabat such that $j'$ increases monotonically with the distance metric from point $\mathbf{d}$. The total reversible heat about this segment is due to the effects as defined in (14) and (15), namely the heat absorbed by the system across the diathermal surface and the heat given up by the temperature well to the reversible heat pumps used to transfer heat to the mass elements. Point $\mathbf{d}$ is connected to $\mathbf{a}$ along a reversible adiabatic pathway and similarly $\mathbf{c}$ to $\mathbf{b}$.

Adiabats $\mathbf{da}$ and $\mathbf{bc}$ are such that we allow $(\mathbf{d} \to \mathbf{a}, \mathbf{c} \to \mathbf{b})$ but the paths $\mathbf{ab}$ and $\mathbf{dc}$ do not coincide in this limit. Along the loop transition $\mathbf{a}-\mathbf{b}-\mathbf{c}-\mathbf{d}-\mathbf{a}$, the net effect would be the summation of all the heat effects in the cycles $S_1, S_2, ....S_N$. We allow for one mass exchange increment of a particular species $\delta m_j$ (for species $j$) for each cycle $S_i$. We define functions $f_{m_j}(i), (j=1,..n; i=1,2,...N')$ such that

$$f_{m_j}(i) = \begin{cases} 1 \text{ if } \delta m_j \text{ is added to the engine } S_i \\ -1 \text{ if } \delta m_j \text{ is removed to the engine } S_i \\ 0 \text{ if no change of } m_j \text{ is noted in engine } S_i. \end{cases}$$

(21)

Then the final value of any component $j$ during the $\mathbf{ab}$ transition at $\mathbf{b}$ will be a finite number where

$$m_j' = m_{j,\mathbf{a}} + \delta m_j \sum_{i=1}^{N'} f_{m_j}(i) \qquad (22)$$

In the limit $\delta m_j \to 0$, all values of $m_j'$ is achievable for arbitrarily large $N'$. Hence the path is completely arbitrary in $\mathbf{m}$ subspace. Then for $N'$ transitions, the result for each of the Carnot cycles $S_i$ means that

$$\sum_{i=1}^{N'}\left(\frac{\delta Q_{tot,i}}{T_1}-\frac{\delta Q_{tot,i'}}{T_2}\right)=0 \qquad (23)$$

Since the **m** subspace transition is completely arbitrary, as is also the isothermal **ab** transition, in the limits

$$N' \to \infty, \delta m_j \to 0, m_j{}'= m_{j,\mathbf{a}}+\lim_{N'\to\infty, \delta m_j \to 0}\delta m_j \sum_{i=1}^{N'}f_{m_j}(i)\;,\;(23)\text{ becomes}$$

$$\oint d\mathcal{S}=\oint \frac{dQ_{tot}}{T}=0. \text{ Hence the state function } \mathcal{S}=\mathcal{S}(\mathbf{P},\mathbf{V},\mathbf{m}) \text{ exists.} \blacksquare$$

In Theorem 1, $dQ_{tot}$ is a composite heat term and from (14,20) we can deduce that

$$d\mathcal{S} = \frac{dQ_{dia}}{T}+\sum_{i=1}^{n}\Delta \mathcal{S}_{i,ss}^{(T_1,p_i)}dm_i \qquad (24)$$

where

$$dQ_{dia} = dQ_{sys}+\sum_{i=1}^{n}dQ_{inj,i} \qquad (25)$$

The total diathermal heat transfer increment by reversible conduction $dQ_{dia}$ consists of $dQ_{sys}$ which refers to the heat absorption by the diathermal boundary of the system when there is no mass transfer and $dQ_{inj,i}$ when there is transfer of substance $j$. The differential of $\mathcal{S}=\mathcal{S}(\mathbf{P},\mathbf{V},\mathbf{m})$ when compared to (24) leads to the following because we have carried out the injection of substances at constant extensive quantities **V**:

$$\left.\frac{\partial \mathcal{S}}{\partial m_i}\right|_{\mathbf{P},\mathbf{V},m_j \neq m_i}=\Delta \mathcal{S}_{i,ss}^{(T,p_i)}, \left.\frac{\partial \mathcal{S}}{\partial V_i}\right|_{\mathbf{P},\mathbf{m},V_j \neq V_i}=\frac{\partial Q_{inj,i}}{\partial V_i},$$

$$\left.\frac{\partial \mathcal{S}}{\partial P_i}\right|_{\mathbf{V},\mathbf{m},P_j \neq P_i}=\frac{\partial Q_{sys}}{\partial P_j}. \qquad (26)$$

In this work, the system is made to absorb heat subjected to extensive variables $\mathbf{V},\mathbf{m}\,(m_i \neq m_j)$ held constant for transfer of mass $m_j$; this concept is different from that provided by Biot [2] since the increment of the heat of mixing/reaction $dQ_{inj,i}$ applies to both chemical reactions and to the heat of mixing, since this heat term arises from the thermalization of intermolecular forces acting on the particles. Biot (ref. 2, p.184, Sec.1) speaks of "a new *intrinsic heat of reaction*, which excludes the heat of mixing, and is more representative of the chemical energy than the standard concept..." In our case, the only distinguishing feature is that in a reversible chemical system which is a Carnot device, the presence of equilibrium constants reduces the number of independent species in **m**, and so conforms in some measure to standard concepts.

(ii) Derivation of a fundamental entropy theorem pertaining only to the system

Theorem 1 is not local in the sense that convected quantities feature in the state function. It would be of singular importance to derive another form which is related only to the system only.

The functions $\Delta \mathcal{S}_{i,ss}^{(T,p_i)}$ refer to the entropy change relative to the standard state ss $(p_0, T_0)$, and $p_i = p_i(\Sigma)$, and are therefore functions of state relative to $\Sigma$. Eq.(24) states

$$\oint \frac{dQ_{dia}}{T} = \oint \sum_{i=1}^{n} \Delta \mathcal{S}_{i,ss}^{(T_1,p_i)} dm_i \qquad (27)$$

It is of interest to examine the r.h.s. circular integral of (27). According to Biot (his eqs. 3.6, 3.17 etc of ref. 2), $ds_T = dq/T$ where $dq$ is the heat absorbed by the primary cell by reversible conduction ( similar to our $dQ_{dia}/T$ ) ".. *is not a state variable* since it does not represent a total entropy change and excludes the part of entropy change due to material convection $dm_k$ into the cell" (ref.2, p188).

The initial state of our elementary system cycle had state $\mathbf{a} = (\mathbf{P'}, T, \mathbf{V}, \mathbf{m})$ ($T$ being the temperature extracted from $\mathbf{P}$). We may map this state to a total convected entropy function $\mathcal{S}_{conv}$ by the function $f_s$, i.e.

$$f_s : (\mathbf{P'}, T, \mathbf{V}, \mathbf{m}) \to \mathcal{S}_{conv} = \sum_{i=1}^{m} \Delta \mathcal{S}_{i,ss}^{(\mathbf{P'},T,\mathbf{V},\mathbf{m})} m_i \qquad (28)$$

where in particular $\mathbf{m}$ is a function of state and the pure state convected entropy functions $\Delta \mathcal{S}_{i,ss}^{\Sigma}$ is a function of state relative to $\Sigma = (\mathbf{P'}, T, \mathbf{V}, \mathbf{m})$, which maps to $(p_i, T)$ of the pure substance state by the function $M_k : \mathbf{a} \to (p_k, T)$ where there is a clear reduction of dimension, since the limits of the integration to yield $\Delta \mathcal{S}_{i,ss}^{(\mathbf{P'},T,\mathbf{V},\mathbf{m})}$ is a function of $\Sigma$ and is not independent therefore of the state $\mathbf{a}$.

**Lemma 5**: The map $f_s$ uniquely identifies the convected entropy for state $\mathbf{a}$.
*Proof*
This may be seen by the method of Gibbs integration[3]. We may physically divide the system $\mathbf{a}$ of our primary cell by its extensive variables by the scale factor $k$ (so that extensive properties $\mathcal{U}$ like the energy of the system or entropy $\mathcal{S}_{conv}$ scale according to $\mathcal{U}_{scale} = k\mathcal{U}(\Sigma = \mathbf{a})$ ) without affecting the intensive function such as $\Delta \mathcal{S}_{i,ss}^{(\{\mathbf{P'},T,\mathbf{V},\mathbf{m}\}=\mathbf{a})}$ because only the intensive variables of state $\mathbf{a}$ and the ratios of the extensive variables determine the state. Hence integrating from $k = 0$ to $k=1$ for state $\mathbf{a}$ for the expression $\int_0^1 \sum_{i=1}^{n} \eta_i k \Delta \mathcal{S}_{i,ss}^{\Sigma} dk$, where $\eta_i = m_i / \sum_{j=1}^{n} m_j$ implies that $\mathcal{S}_{conv}$ is the total convected entropy for $\mathbf{a}$ relative to the initial state of $\mathbf{0}$ extension whose entropy is set at zero. Suppose $\mathcal{S}_{conv}$ is not unique for the state, so that (a) for the same state there are two different values $\mathcal{S}_{conv}(1)$ and $\mathcal{S}_{conv}(2)$, or (b) that the intensive variables $\Delta \mathcal{S}_{i,ss}^{(\{\mathbf{P'},T,\mathbf{V},\mathbf{m}\}=\mathbf{a})}$ are different for the same value of $\mathcal{S}_{conv}$, then $\mathcal{S}_{conv}(1) - \mathcal{S}_{conv}(2) = \sum_{i=1}^{n}(\Delta\mathcal{S}_{i,ss}^{\Sigma}(1) - \Delta\mathcal{S}_{i,ss}^{\Sigma}(2))m_i$, and this expression is zero since the specific entropies are functions of state at $\Sigma$, contradicting the hypothesis, so that the mapping is unique ■

**Theorem 2**:
The non independence of the functions $\Delta\mathcal{S}_{i,ss}^{(\{\mathbf{P'},T,\mathbf{V},\mathbf{m}\})}$ is given by the differential expression $\sum_{i=1}^{n} d(\Delta\mathcal{S}_{i,ss}^{(\{\mathbf{P'},T,\mathbf{V},\mathbf{m}\})})m_i = 0$ and furthermore, $d\mathcal{S}_{conv} = \sum_{i=1}^{m} \Delta\mathcal{S}_{i,ss}^{\mathbf{a}} dm_i$ is the differential of the state function $d\mathcal{S}_{conv}$.

*Proof*
Since $m_i$ and $\Delta\mathcal{S}_{i,ss}^{\Sigma}$ are functions of state, then so is $\mathcal{S}_{conv}$. The increment of convected entropy during a system transition has been defined (7) as
$dS_{conv} = \sum_{i=1}^{n} \Delta\mathcal{S}_{i,ss}^{\Sigma} dm_i$.    On    the    other    hand    (28)    yields
$dS_{conv} = \sum_{i=1}^{n}\left(m_i d\Delta\mathcal{S}_{i,ss}^{\mathbf{a}} + \Delta\mathcal{S}_{i,ss}^{\mathbf{a}} dm_i\right)$. Comparing the preceding two expressions for $d\mathcal{S}_{conv}$ proves the theorem ■

The first result of Theorem 2 is very similar in form to the Gibbs-Duhem type equations (ref.3, Chapter 2, p.92-106).

In view of Theorem 2, we can state Theorem 3.

**Theorem 3:**
There exists an entropic function of state $\mathcal{S}_{dia}$ with differential given by $dS_{dia} = \frac{dQ_{dia}}{T}$ which is distinct from the total system entropic function of Theorem 1 where $Q_{dia}$ is a function representing the total heat absorption of the system through a diathermal boundary. In particular, this state function is local in the sense that it refers to quantities measured at the system during system transitions.
*Proof*
From (27), since the r.h.s. integral is zero from Theorem 2, the result follows ■

From Theorem 3, we infer that there are two calorimetrically distinct entropies for open systems, and that the Gibbs' method for open systems does not distinguish these, leading to ambiguous if not incorrect assignments of heat transfer terms. In particular, Theorem 3 contradicts Biot's assertion as mentioned earlier, and it represents a new fundamental result in open-system thermodynamics.

## 5. Discussion of Traditional Open-System Thermodynamics
For the energy $\mathcal{U}$, its differential in traditional thermodynamics is written
$d\mathcal{U} = Td\mathcal{S} - pdV + \sum_{i=1}^{N} \mu_i dn_i$     (29)

where Denbigh explicitly states (ref.3, p.81-82) that the above refer only to closed systems; for open systems the term $Td\mathcal{S}$ "..is no longer interpretable as heat". On the other hand, (ref.3, p.78,eq. 2.40), in defining the chemical potential (at constant entropy), it is clear that for Denbigh's rationalization -based on Gibbs' definition- to be consistent, $d\mathcal{S}$ must refer to the *total* entropy increment which does not involve merely heat transfer across a diathermal boundary but the other terms of Theorem 1. However, we show here that $d\mathcal{S}$ arises from the sum of diathermal heat transfer and another convective heat term $\Delta Q_{total}^{\Sigma}$ of Lemma 1, so when appropriately defined, $Td\mathcal{S}$ can refer to a generalized heat term. The entire field of irreversible thermodynamics and simulations have viewed the terms as in (29) as pertaining to a local system without the external convective terms; this is also true of purported generalizations of the basic open system equations in "extended irreversible thermodynamics" where this locality is assumed[4]; in particular the boundary is a key element in demarcating the system from the environment or supply elements but these distinctions are usually not considered in the local equilibrium hypothesis used for non-equilibrium extensions[5]. Quite remarkably, perhaps, equilibrium specific enthalpy data from NEMD simulations based on Onsager reciprocity have been calculated with a heat flux term defined without recourse to system boundaries , and based on the presumed local nature of these open system equations. Such results may imply a remarkable cancellation of terms and the relative insignificance of external convective factors in the final computations[6] if these results are correct.

Here, the local energy differential $d\mathcal{U}$ at point **a** of the system $S$ only is

$$d\mathcal{U}(\mathbf{a}) = Td\mathcal{S}_{dia} - \mathbf{P}.\mathbf{dV} + \sum_{i=1}^{n} \mu_i dm_i \qquad (30)$$

where the "chemical potential" $\mu_i$ is a mixture of heat and work terms

$$\mu_i = \Delta W_i^{form} + W_{ss,i} + \int_{ss}^{\mathbf{a}} T' d\Delta \mathcal{S}_{i,ss}^{(T',p_k')} + \zeta_i(\mathbf{a}) + \mathcal{U}_{ss,i} .$$

$\mathcal{U}_{ss,i}$ is the standard state energy assignment for the pure substance. It will be noticed that the total entropy $\mathcal{S}$ has a convected entropy term missing in many standard treatments, and in particular the chemical potential here has got a thermal energy component often missing in treatments that write $\mu_i$ in the form $\mu_i = C^{\ominus} + RT\ln a_i$, where clearly only isothermal work terms feature since $C^{\ominus}$ is almost always considered a function of $T$ only, and only in ideal systems would the temperature reflect the thermal energy, independent of concentration of the species. The thermal energy changes could conceivably be absorbed by the activity coefficient term but its definition in statistical mechanics given below seems to preclude this possibility. Eq.(30) does not have any "convected" entropy entering into the local energy equation; in particular, $d\mathcal{U}(\mathbf{a})$ is dependent of the energy zero of the pure substances, $\mathcal{U}_{ss,i}$. Standard statistical mechanics define the activity coefficient $\gamma_i$ of species $i$ as the $i$'s concentration or density derivative of the excess Helmholtz function $A''$ [7](or in other instances the excess Gibbs free energy function) when intermolecular forces exist with the ideal state as the reference, where $a_i = \gamma_i c_i$ with $c_i$ the density or concentration of species $i$. All these standard generalized definitions may not fully account for the thermal energy increments and seem to account for the average work contribution only at fixed pressure or volume element, such as what is obtained in elementary Debye-Huckel theory of activity coefficients and its manifold extentions[8], because of the form of the energy expression, where all heat transfer is delegated to the $Td\mathcal{S}$ term, which does not refer to material transfer if one interprets the term as a reversible thermal energy transfer. However more examination is required before any definitive conclusions may be drawn. Where entropy is concerned, the overall entropy of Theorem 1 is non-local, whereas that of Theorem 3 is; when some central, non-equilibrium local entropy is postulated to have certain properties without proof in new theories (see ref.4, pg. 42) it is difficult to understand what is being meant, and whether there is an existence theorem to guarantee that such an entropy does in fact exist.

## 6. Discussion of Internal and External Processes

In view of Theorem 3, it is possible to separate in some instances processes external to the system or cell which is enclosed by a diathermal boundary, and the processes within the boundaries of the system. It is of interest to formulate their relationship.

(i) *Relationship between a closed Carnot cycle and an Open system working in a cycle*

Consider a two-port mass exchange (either isothermal or adiabatic) Carnot cycle j that is infinitely small in extent and heat delivery and whose endpoints lie on the isothermal paths used to prove Theorem 1. Let this system absorb heat $Q_1^j$ at temperature $T_2$ and give up heat $Q_2^j$ at $T_1$ in a cycle at the isothermal segments. Let mass of substance $i$ introduced at the input port be $\delta m_{j,i}$ having total energy $\mathcal{U}_{in,j,i}$ before input and let the mass have energy $\mathcal{U}_{out,j,i}$ on output. This elementary cycle in the limit belongs to a closed loop where $\oint d\mathcal{S}_{dia} = \oint \frac{dQ_{dia}}{T} = 0$. For cycle $j$, the relations where the system is concerned is obviously

$$\Delta \mathcal{U}_{j,i} + Q_1^j - Q_2^j - W_{ext}^j = 0 \qquad (31)$$

where $\frac{Q_1^j}{Q_2^j} = \frac{T_1}{T_2}$ and $\Delta \mathcal{U}_j = \sum_{i=1}^{n} \mathcal{U}_{in,j,i} - \mathcal{U}_{out,j,i}$. $W_{ext}^j$ is the total work done by the system, including the work of injection/extraction. Only in a closed system cycle would $\Delta \mathcal{U}_j = 0$ and the virtual work done by the system $W_{vir}^j$ is defined as

$$W_{vir}^j = Q_1^j \left( \frac{T_1^j - T_2^j}{T_1^j} \right). \qquad (32)$$

Comparing (31) and (32) yields

$$Q_1^j - Q_2^j = W_{vir}^j = W_{ext}^j - \Delta \mathcal{U}_j. \qquad (33)$$

Summing back (33) as $N \to \infty$, we get the defining relations for increments

$$\oint d\mathcal{S}_{dia} = \sum_{j=1}^{N}\frac{Q_1^j}{T_1} - \frac{Q_2^j}{T_2} = 0, \sum_{j=1}^{N} Q_1^j - Q_2^j = \oint dQ_{dia} \quad , \qquad \oint dW_{vir} = \sum_{j=1}^{N} W_{vir}^j \quad \text{and}$$

$$\oint d\mathcal{U}_i = \sum_{j=1}^{N}\left(\mathcal{U}_{in,j,i} - \mathcal{U}_{out,j,i}\right),$$ where in particular $d\mathcal{U}_i$ refers to the energy of species $i$ introduced into the system. So the relation in the limit becomes

$$\oint dW_{vir} = \oint dW_{ext} - \oint \sum_{i=1}^{n} d\mathcal{U}_i \qquad (34)$$
$$= \oint dW_{ext}$$

since $\oint \sum_{i=1}^{n} d\mathcal{U}_i = 0$. Result (34) may be summarized in terms of a **Principle of correlation between a closed and open Carnot cycle:**
*In an open system exchanging heat with reservoirs, the total external work performed by the system on the environment only in a cycle is equivalent to that performed by a closed Carnot engine exchanging the same amount of heat with the reservoirs along a cycle.*

(ii) *The nonequivalence of the total system-environmental entropy $d\mathcal{S}$ and the local entropy $d\mathcal{S}_{adia}$ differentials*

On completing a loop $a \xrightarrow{P1} b \xrightarrow{P2} a$ where paths $P1$ and $P2$ are different, and where $a, b$ represent points in the thermodynamic phase space, in the case of the variables used to compute $\mathcal{S}$, they are all brought back to the standard states in $Cp_i$ in any elementary cycle but this is not the case with the variables used for computing $\mathcal{S}_{adia}$, where in a cycle, the environmental mass elements from $Cp_i$ need not return back to the standard states. One might therefore predict different global properties arising from considerations of these two types of entropies. At any point during the transition, $d\mathcal{S}_{adia} = \frac{dQ_{adia}}{T}$, but for the composite system, we have from (24) that $d\mathcal{S} = \frac{dQ_{dia}}{T} + \sum_{i=1}^{n} \Delta \mathcal{S}_{i,ss}^{(T_1, p_i)} dm_i$, so that $d\mathcal{S} \neq d\mathcal{S}_{adia}$. One might therefore state a **Principle of open Carnot cycle heat conversion**:
*The total heat absorbed in a reversible arbitrary cycle of an open Carnot engine must equal the total work done by the engine on the environment and the energy change associated with the movement of convected substances outside the system, where the diathermal entropy change of the system is zero, but not necessarily of the environment.*

## 7. Conclusions
We have shown from first principles that various new functions of states exists, especially in entropy, that may be applied rigorously to open systems. They were all derived analytically without presuming the existence of functions of state, which is the normal practice, even if it is not analytical. The present formulation when contrasted with Gibbs' does not suffer from ambiguity of definition and circular reasoning which some writers like Biot claim is evident from the Gibbsian approach.

**acknowledgements**: I am grateful for the over three decades of solitary and seminal work that Biot put into his thermodynamics of hypersystems; I have profited from his definitions and structures of thought, although the results here have little to do with his.